\documentclass{article}
\usepackage{amssymb,latexsym,amsmath,amsthm}
\usepackage{fancyhdr}

\usepackage{color}
\usepackage{mathrsfs}

\usepackage{lipsum}

\newcommand\blfootnote[1]{%
  \begingroup
  \renewcommand\thefootnote{}\footnote{#1}%
  \addtocounter{footnote}{-1}%
  \endgroup
}

\numberwithin{equation}{section}
\newtheorem{lemma}{Lemma}
\newtheorem{theorem}[lemma]{Theorem}
\newtheorem{corollary}[lemma]{Corollary}
\newcommand{\ep}{\varepsilon}

\newcommand{\Z}{\mathbb{Z}}
\newcommand{\Q}{\mathbb{Q}}

\newcommand{\beql}[1]{\begin{equation}\label{#1}}
\newcommand{\eeq}{\end{equation}}

\newcommand{\modd}[1]{\; ( \text{mod} \; #1)}

\newcommand{\Cl}{\mathrm{Cl}}

\newcommand{\card}{\#}

\begin{document}

\title{$\ell$-Torsion in Class Groups via Dirichlet $L$-functions}  
\author{D. R. Heath-Brown\\
Mathematical Institute, Oxford }
            
  \date{}
  \maketitle
  
\begin{abstract}
For a prime $\ell$, let $h_\ell(K)$ denote the $\ell$-part of
the class number of the number field $K$. We investigate upper bounds
for $h_\ell(K)$ when $K$ is
quadratic or cubic, particularly in the case in which the
discriminant of $K$ is smooth. This is achieved using properties
of Dirichlet $L$-functions.
\end{abstract}
\noindent
\newpage

\section{Introduction}
Let $K$ be an algebraic number field of degree $n$, and write
$\Delta_K$ for its discriminant, $\Cl(K)$
for its class group, and $h(K)=\card\Cl(K)$
for its class number. It then follows from a theorem of
Landau\footnote{The weaker bound
  $h(K)\ll_{n,\ep}|\Delta_K|^{1/2+\ep}$ is sometimes described as
  coming from the Brauer--Siegel Theorem. Unfortunately the size of
  this footnote is not sufficient to list all the reasons why this
  mis-attribution is wrong.}
\cite{Landau} from 1918 that
\[h(K)\ll_n |\Delta_K|^{1/2}(\log|\Delta_K|)^{n-1},\]
where $\Delta_K$ is the discriminant of $K$.
For a prime $\ell$ we write $\Cl_\ell(K)$ for the
$\ell$-torsion subgroup of $\Cl(K)$, that is to say the subgroup
consisting of ideal classes whose order divides $\ell$. Writing
$h_\ell(K)=\card\Cl_\ell(K)$, Landau's theorem produces the ``trivial bound''
\[h_\ell(K)\ll_n |\Delta_K|^{1/2}(\log|\Delta_K|)^{n-1}.\]
However it is natural to conjecture a stronger bound, taking the form
\beql{h-conj}
h_\ell(K) \ll_{n,\ell,\ep} |\Delta_K|^{\ep} 
\eeq
for any $\ep>0$. One can prove this when there is an appropriate genus
theory, see Kl\"{u}ners and Wang \cite{KW},
and in particular  when $K$ is a Galois extension of $\Q$ of degree $\ell$.
Thus when $K$ is quadratic, for example, one has the classical estimate
\[h_2(K)\ll 2^{\omega(|\Delta_K|)}\ll_\ep|\Delta_K|^{\ep} ,\]
where $\omega(d)$ counts distinct prime factors of $d$.

Aside from cases where genus theory applies in some form or other, we
cannot yet prove the 
conjecture (\ref{h-conj}), and so we may ask whether there is any
nontrivial bound of the shape
\[h_\ell(K) \ll_{n,\ell} |\Delta_K|^{1/2-\delta},\]
for some $\delta>0$. There are various results that should be mentioned
here. Firstly, Ellenberg and Venkatesh \cite[Corollary 3.7]{EV},
have shown that
\[h_3(K)\ll_\ep|\Delta_K|^{1/3+\ep} \]
for any $\ep>0$, whenever $K$ has degree 2 or 3, and that there is
some $\delta>0$ such that
\[h_3(K)\ll |\Delta_K|^{1/2-\delta} \]
whenever $K$ has degree 4. Secondly, Bhargava, Shankar, Taniguchi,
Thorne, Tsimerman and Zhao \cite{BSTTTZ} have shown that
\[h_2(K) \ll_{n,\ep} |\Delta_K|^{1/2-\delta_n+\ep},\]
with $\delta_n=0.2215\ldots$ for $n=3$ or 4, and with $\delta_n=1/2n$
for $n>4$. Finally, Wang \cite{JW1} has shown that, for any primes $\ell$ 
and $p$, there is a corresponding $\delta(\ell,p)>0$ such that
\beql{WB}
h_\ell(K)\ll |\Delta_K|^{1/2-\delta(\ell,p)} 
\eeq
whenever $K/\Q$ is Galois, with Galois group $(\Z/p\Z)^r$ for some $r\ge 2$.
Indeed her result is somewhat more general than this, and in a later preprint
\cite{JW2} more general Galois groups are treated. 

One can do more under the Generalized Riemann Hypothesis (GRH).
With this assumption
Ellenberg and Venkatesh \cite[Proposition 3.1]{EV}, show that
\beql{EVB}
h_\ell(K) \ll_{n,\ell,\ep} |\Delta_K|^{1/2-1/2\ell(n-1)+\ep},
\eeq
for any $\ep>0$, and any prime $\ell$. (Note however that
with $n=p^r$, Wang's bound (\ref{WB}) has $\delta(\ell,p)>1/2\ell(n-1)$
for large enough $r$.)
Our goal in this paper is to
investigate the extent to which one can relax, or indeed remove, the
hypothesis of the GRH from (\ref{EVB}), particularly in the case of 
quadratic and cubic fields.

We begin by mentioning one almost trivial result. This considers
``pure" cubic fields $K=\Q(\sqrt[3]{d})$, where $d>1$ is a
cube-free integer.
\begin{theorem}\label{PureThm}
  For a pure cubic field we have
\[h_\ell(K) \ll_{\ell,\ep} |\Delta_K|^{1/2-1/(4\ell)+\ep},\]
for any $\ep>0$, and any prime $\ell$.
\end{theorem}
The principle behind this bound has been known for some time,
in lectures by Wang in particular, but we feel the result merits explicit
mention here.

In the context of Theorem \ref{PureThm} we should note that in the
case $\ell=3$ one has the stronger bound
\begin{equation}\label{ell3}
  h_3(K) \ll 9^{\omega(d)}\ll_{\ep}|\Delta_K|^{\ep},
\end{equation}
for any $\ep>0$ and any pure cubic field $K$,
where $\omega(k)$ is the
number of distinct prime factors of $k$. This follows from Gerth
\cite[Theorem 2.7 \& Proposition 3.2]{Gerth}. In the notation of
\cite{Gerth}, the 3-rank of $\Cl(K)$ is denoted $\mathrm{rank}\; S_L$, and
according to Theorem 2.7 this is at most a certain quantity
$\mathrm{rank}\;{}_ND$. Proposition 3.2 gives a formula for 
$\mathrm{rank}\;{}_ND$, which shows that it is at most $r+t-1$, where
$r$ is the 3-rank of the class group of $\Q(\sqrt{-3})$ in our case
(so that $r=0$), and $t$ is the number of ramified primes in the
extension $\Q(\sqrt[3]{d},\sqrt{-3})/\Q(\sqrt{-3})$. Thus $t$ is at
most $2\omega(d)+1$, and (\ref{ell3}) follows.

Ellenberg and Venkatesh's proof of (\ref{EVB}) builds on ideas of
Soundararajan \cite{sound} to give the following result, which is a special case
of Ellenberg and Venkatesh \cite[Lemma 2.3]{EV}.
\begin{lemma}\label{FundLem}
Let $K$ be a number field of degree $n$, let $\ep>0$, let $\ell$ be a prime, and
let $\varpi<1/2\ell(n-1)$. Suppose there are $M$ unramified prime
ideals of $K$, all of degree 1, and having norm at
most $|\Delta_K|^\varpi$. Then
\[h_\ell(K) \ll_{n,\ell,\ep,\varpi} |\Delta_K|^{1/2+\ep}M^{-1}.\]
\end{lemma}

Soundararajan restricted his attention to imaginary quadratic
fields. The extension to general number fields by Ellenberg and
Venkatesh required a significant new idea to handle infinite unit
groups. 

To deduce Theorem \ref{PureThm} from Lemma \ref{FundLem}
it suffices to observe that if
$p\equiv 2\modd{3}$ is any prime not dividing $\Delta_K$, then $(p)$
splits over $K$ as a product $P_1P_2$ of distinct prime ideals, where
$P_1$ has degree 1. This provides $M\gg|\Delta_K|^\varpi/\log|\Delta_K|$
suitable prime ideals, and the result follows on taking
$\varpi=1/(4\ell)-\ep$, and re-defining $\ep$.
\bigskip

In order to obtain (\ref{EVB}) from Lemma \ref{FundLem} Ellenberg and
Venkatesh call on a strong form of the Chebotarev Density Theorem
which assumes the GRH. The principal theme of this paper is to ask
whether one can get by without assuming such a strong hypothesis
as GRH. We shall confine our attention to situations in which it
suffices to work with Dirichlet $L$-functions, although some of the ideas
apply to more general Artin $L$-functions. This means that we will
assume $K$ to be either quadratic or cubic. 

When $K$ is quadratic we
might seek rational primes $p$ for which $\chi(p)=+1$, where
$\chi$ is the quadratic character associated to $K$.
Any such prime splits in $K$ into two unramified first-degree prime
ideals.  The study of such primes is hampered by the possibility of 
$L(s,\chi)$ having an exceptional real zero, 
very close to $s=1$. We can largely avoid this via the following easy
variant of Lemma \ref{FundLem}.
\begin{lemma}\label{FundLemV}
Let $K$ be a number field of degree $n$, let $\ell$ be a prime, let $\ep>0$ and
let $\varpi<1/2\ell(n-1)$. Suppose there are $M$ integral
ideals $I$ of $K$,  having square-free norm $N(I)$ coprime to $\Delta_K$, and
with $N(I)\le|\Delta_K|^\varpi$. Then
\[h_\ell(K) \ll_{n,\ell,\ep,\varpi} |\Delta_K|^{1/2+\ep}M^{-1}.\]
\end{lemma}
We will prove this later, in Section \ref{QF}.  Lemma \ref{FundLemV} allows 
us to work with $L(s,\chi)$ directly, rather than with its logarithmic
derivative.
We will prove the following result.
\begin{theorem}\label{QT}
Let $\ell\ge 3$ be prime. Let $\theta\in(0,\tfrac12)$ and
\[0<\xi<\frac{2\theta-\theta^2}{4\ell}.\]
Suppose that $K$ is a quadratic number field with corresponding quadratic
character $\chi$, and assume that
\[L(\sigma+it,\chi)\ll |\Delta_K|^{\xi},\;\;\; (\sigma\ge 1-\theta,\;
|t|\le 1).\]
Then 
\[h_\ell(K) \ll_{\ell,\theta,\xi,\ep} |\Delta_K|^{1/2-1/(2\ell)+\ep},\]
for any fixed $\ep>0$. The implied constant above is ineffective.
\end{theorem}
The Generalized Lindel\"{o}f Hypothesis would permit us to
take any $\xi>0$, for any $\theta$, so that the assumptions required
above are vastly weaker even that this. Indeed an examination of the proof
shows that it would suffice to have a suitable bound on $L(s,\chi)$ on a rather 
shorter interval than $|t|\le 1$. We also remark that the proof can be adapted
to fields $K$ of higher degree, given a suitable assumption on the size of
$\zeta_K(s)/\zeta(s)$.
\bigskip

When $K$ is a cyclic cubic field the way $p$ splits is
determined by a cubic character $\chi_K$ whose conductor is composed
of primes that divide $\Delta_K$. In this situation any prime $p$
for which $\chi_K(p)=1$ will split over $K$ into three unramified
first-degree prime ideals. (The character $\chi$ is only defined
up to conjugation.) On the other hand, for a non-cyclic cubic field the
discriminant $\Delta_K$ will be a non-square, and any prime $p\ge 3$
with
\[\left(\frac{\Delta_K}{p}\right)=-1\]
factors over $K$ as a product of two prime ideals, one having degree 1
and the other with degree 2. Thus it is enough for us to find a good
supply of rational primes $p$ for which $\chi(p)$ has a prescribed
value, where $\chi$ is a quadratic or cubic character. The conductor $q$
of $\chi$ only involves primes dividing $\Delta_K$, and the argument
of Heath-Brown \cite[foot of page 271]{HBlinnik} allows us to deduce that
$q\ll|\Delta_K|$. Our situation is similar to Vinogradov's problem on
the least quadratic non-residue modulo $q$, but our task involves
finding many suitable primes, not just one. We attack this by
considering zero-free regions of $L$-functions.

\begin{theorem}\label{CT}
Let $\ell$ be prime, and let $\ep>0$. Then there is an effective constant
$C_1(\ell,\ep)$ with the following property.

Let $K$ be a cubic number field. If $K$ is 
cyclic let $\chi$ one of the primitive cubic characters associated to $K$, and 
if $K$ is non-cyclic cubic field let $\chi$ be the primitive quadratic 
character which induces the Jacobi symbol $(\Delta_K/*)$.

Suppose that $1\le\vartheta\le\log|\Delta_K|$ and that $L(s,\chi)$ 
has no zeros in the disc $|s-1|\le\vartheta/\log|\Delta_K|$, apart possibly 
from a single real zero in the case in which $\chi$ is quadratic. Then if 
$\vartheta\ge C_1(\ell,\ep)$ one has
\[h_\ell(K) \ll_{\ell,\vartheta,\ep} |\Delta_K|^{1/2-1/(4\ell)+\ep}\]
with an effective implied constant.
\end{theorem}
We should emphasize that the assumption here is vastly weaker 
than the Generalized Riemann Hypothesis. The standard zero-free
region produces an absolute constant $C>0$ such that the
disc $|s-1|\le C/\log q$ has the required property, where $q$ is the
conductor of $\chi$. Hence, for example, we may generalize
Theorem \ref{PureThm} to cover any sequence of non-cyclic cubic
fields $K$ whose discriminants factor as $\Delta_K=a_Kb_K^2$ with
$\log|a_K|/\log|\Delta_K|\to 0$.
\bigskip

We may obtain zero-free regions via a $q$-analogue of a result which 
is familiar in the context of the Riemann Zeta-function, see Titchmarsh
\cite[Theorem 3.10]{Titch}, for example.
\begin{theorem}\label{zfr}
There is an effectively computable constant $C_2>0$ with the following
property.  Suppose that $\chi$ is a complex Dirichlet character
modulo $q$, and that the estimate
  \beql{zfrc}
  |L(\sigma+it,\psi)|\le e^{\phi},\;\;\;
  (\forall\sigma\ge 1-\theta,\;\;\forall |t|\le 3),
\eeq
holds both for $\psi=\chi$ and for $\psi=\chi^2$, where $\theta$ and
$\phi$ are parameters satisfying $\phi e^{-\phi}\le\theta\le 1\le\phi$
and $\phi\ge\theta(1+\log\log 3q)$.
Then $L(s,\chi)$ has no zeros with
\beql{ZFR}
\sigma\ge 1-C_2\theta/\phi,\;\; |t|\le 1.
\eeq

Moreover, if $\chi$ is a character modulo $q$ of order 2, and the
bound (\ref{zfrc}) holds for $\psi=\chi$, then $L(s,\chi)$ has at most
one zero in the region (\ref{ZFR}). 
If this zero exists, it is a simple real zero of $L(s,\chi)$.

Finally, if $\chi_1$ and $\chi_2$ are distinct characters modulo $q$,
both of order 2, such that the
bound (\ref{zfrc}) holds for $\psi=\chi_1$, $\psi=\chi_2$ and
$\psi=\chi_1\chi_2$, then $L(s,\chi_1)$ and $L(s,\chi_2)$ cannot both have
zeros in the region (\ref{ZFR}). 
\end{theorem}
The proof is entirely straightforward and we do not claim any
novelty. See Graham and Ringrose \cite[Lemma 6.2]{GR} for a special
case of the result.

As an immediate corollary of Theorems \ref{CT} and \ref{zfr} we have 
the following.
\begin{corollary}\label{cor}
For each prime $\ell\ge 3$ and small real number $\ep>0$ there is an 
effective constant $C_3(\ell,\ep)$ with the following property. 

Let $K$  and $\chi$ be as in Theorem \ref{CT}, and suppose that
\[L(\sigma+it,\chi)\ll |\Delta_K|^{\xi},\;\;\; (\sigma\ge 1-\theta,\;
|t|\le 3)\]
for some $\theta\in(0,1)$ and $\xi\le C_3(\ell,\ep)\theta$. Then
\[h_\ell(K) \ll_{\ell,\theta,\xi,\ep} |\Delta_K|^{1/2-1/(2\ell)+\ep}.\]
The implied constant above is effective.
\end{corollary}
For the proof it is enough to take $C_3(\ell,\ep)=C_2/C_1(\ell,\ep)$, 
which will be sufficient when $|\Delta_K|$ is large enough.

It is interesting to compare Theorem \ref{QT} with Corollary \ref{cor}.
The former is ineffective, but has a precise condition on $\xi$, while 
the latter is effective, with an inexplicit requirement for $\xi$. The 
distinction between the ranges $|t|\le 3$ and $|t|\le 1$ has no 
significance. One could easily prove the results with either range.
\bigskip

One situation in which we are able to say something useful about
bounds for $L(s,\chi)$ is that in which we restrict to smooth values of $q$.
The following is a trivial variant of Theorem 5 of Graham and Ringrose 
\cite{GR}.
\begin{theorem}\label{CharSum}
Let $\chi$ be a non-principal character to modulus $q$. Let $p$ be
the largest prime factor of $q$ and write $r$ for the largest 
square-full factor of $q$. Let $k\in\mathbb{N}$ and write $L=2^{k+3}-2$. Then
\begin{eqnarray*}
\lefteqn{\sum_{M<n\le M+N}\chi(n)}\\
&\ll_r & M^{1-(k+3)/L}
  q^{1/L}d(q)^{(3k^2+11k+8)/(2L)}(\log q)^{(k+3)/L}\sigma_{-1}(q)p^{(k^2+3k+4)/(4L)}
  \end{eqnarray*}
  for any positive integers $N\le M$.
  The implied constant above depends only on $r$, and is effective.
\end{theorem}

This is proved by the $q$-analogue of van der
Corput's method. It corresponds to the $(k+3)$-rd derivative
estimate, see Titchmarsh \cite[Theorem 5.13]{Titch} for example.

The procedure for deriving estimates for $L(s,\chi)$ from bounds for
the corresponding character sum is standard, and we immediately get
the following corollary. For our purposes a
relatively crude result is sufficient.
\begin{corollary}\label{Lbound}
Under the hypotheses of Theorem \ref{CharSum}, for any fixed
$\ep>0$ we will have
\[L(s,\chi)\ll_{r,k,\ep}q^{1/L+\ep}p^{(k^2+3k+4)/(4L)}\]
throughout the region $\sigma\ge 1-(k+3)/L$, $|t|\le 3$.
\end{corollary}

We may feed this into Theorem \ref{QT} to
produce the following.
\begin{theorem}\label{HLGR}
  Suppose we are given a prime $\ell$. Then there is a corresponding
  effective constant $\delta_2(\ell)>0$ with the following property. 
  
  Let $K$ be a quadratic number field such that 
  $p\le|\Delta_K|^{\delta_2(\ell)}$ for every
  prime divisor $p$ of $\Delta_K$. Then 
\[h_\ell(K) \ll_{\ell,\ep} |\Delta_K|^{1/2-1/(2\ell)+\ep},\]
for any $\ep>0$, with an ineffective implied constant. For example,
one may take $\delta_2(5)=1/343$.
\end{theorem}

As in Heath-Brown \cite[foot of page 271]{HBlinnik}, if $\chi$ is a character 
of bounded order, its conductor $q$ has bounded square-full part, so that
$r\ll 1$ in the notation of Theorem \ref{CharSum}. If $p\le|\Delta_K|^\delta$ 
then Corollary \ref{cor} allows us to take
\[\theta=\frac{k+3}{2^{k+3}-2}\;\;\;\mbox{and}\;\;\;
\xi=\frac{1+\delta(k^2+3k+4)/4}{2^{k+3}-2}-\ep\]
for any $\ep>0$. One then sees that $\xi/\theta$ can be made arbitrarily 
small, by choosing $k$ large enough and $\delta$ small enough.
Indeed one may check that the values $k=15$
and $\delta=1/343$ are sufficient.

Similarly, combining Corollary \ref{cor} with Theorem \ref{CT}
yields:-
\begin{theorem}
  Suppose we are given a prime $\ell$. Then there is a corresponding
  effective constant $\delta_3(\ell)>0$ with the following property. 
  
  Let $K$ be a cubic number field such that 
  $p\le|\Delta_K|^{\delta_3(\ell)}$ for every
  prime divisor $p$ of $\Delta_K$. Then 
\[h_\ell(K) \ll_{\ell,\ep} |\Delta_K|^{1/2-1/(2\ell)+\ep},\]
for any $\ep>0$, with an effective implied constant. 
\end{theorem}
\bigskip

{\em Acknowledgements.} The author is grateful to Lillian Pierce and
Jiuya Wang for helpful discussions, and in particular for
providing references to relevant works.

\section{Quadratic Fields --- Proof of Theorem \ref{QT}}\label{QF}
We begin by proving Lemma \ref{FundLemV}.
This is done by the obvious modification of the argument given
by Ellenberg and Venkatesh in proving \cite[Lemma 2.3]{EV}. All that is 
necessary for Lemma \ref{FundLemV} is to note that if $I_1$ and $I_2$ are 
ideals as in the lemma, such that $I_1^\ell I_2^{-\ell}=(u)$ for some 
$u$ in a proper subfield $K_1$ of $K$, then $I_1=I_2$. For if $u$ is not a 
unit, then $v_P(u)\not=1$ for some prime ideal $P$ of $K_1$, so that
$P$ divides 
at least one of $I_1$ or $I_2$. This however is impossible since 
\[N_{K/\mathbb{Q}}(P)=N_{K_1/\mathbb{Q}}(P)^{[K:K_1]}\]
whereas $N(I_1)$ and $N(I_2)$ are square-free.
\bigskip

Now let $K$ be a quadratic field, as described in Theorem \ref{QT},
and take $\chi$ 
to be the associated real character, so that $\chi$ has conductor
$q=|\Delta_K|$. 
The parameters $\ell$, $\theta$ and $\xi$ will be fixed throughout
our argument, and 
we allow all our order constants to depend on them without further mention.

We define
\[f(s)=\sum_{n=1}^{\infty}\mu^2(n)\prod_{p|n}(\chi_0(n)+\chi(n))=\
\sum_{n=1}^\infty a_n n^{-s},\]
so that $a_n$ is the number of integral ideals of $K$ of norm $n$, if $n$ is
square-free and co-prime to $\Delta_K$, and $a_n=0$ otherwise. Using 
Euler products we see that we may write
\[f(s)=L(s,\chi_0)L(s,\chi)f_0(s)\]
where $f_0(s)$ has an Euler product which is absolutely convergent for
$\sigma>\tfrac12$, with $f_0(s)\ll 1$ uniformly for $\sigma\ge\tfrac34$, say. 

Our argument now introduces some rather severe Gaussian weights, which 
are needed because we only have information for the narrow range $|t|\le 1$.
We consider
\beql{I2}
I=\frac{1}{2\pi i}\int_{2-i\infty}^{2+i\infty}f(s)\exp( s^2y)ds=
\frac{1}{2\sqrt{\pi y}}\sum_1^\infty a_n \exp\{-(\log n)^2/4y\}
\eeq
with $y=\kappa\log q$, where $\kappa\in(0,1)$ is a constant to be chosen
later. We let $\delta\in(0,1)$ be another constant, also to be chosen later,
 and move the path of
integration so that it runs from $1+\delta-i\infty$ to $1+\delta-i$
to $1-\theta-i$ to 
$1-\theta+i$ to $1+\delta+i$ to $1+\delta+i\infty$. There is a pole at $s=1$, 
which contributes
\beql{pb}
\prod_{p|q}(1-p^{-1})L(1,\chi)f_0(1)e^y\gg_\delta q^{-\delta}e^y,
\eeq
by Siegel's Theorem. The implied constant here 
is ineffective. The section of the integral from $1+\delta+i$ to
$1+\delta+i\infty$ produces
\[\ll_\delta \int_1^\infty \exp((1+\delta)^2y-t^2y)dt\le
e^{(1+\delta)^2y}\int_1^\infty\exp(-ty)dt=
\frac{e^{(2\delta+\delta^2)y}}{y},\]
and similarly for the section from $1+\delta-i\infty$ to
$1+\delta-i$. When 
$1-\theta\le\sigma\le 1+\delta$ and $t=\pm 1$ we have
\[f(s)\ll |L(s,\chi_0) |q^\xi\ll_\delta q^{\xi+\delta}\]
and
\[|\exp(s^2 y)|= \exp(\sigma^2 y-y)\le \exp((2\delta+\delta^2)y),\]
whence the integrals from $1+\delta-i$ to $1-\theta-i$ and from 
$1-\theta+i$ to $1+\delta+i$ contribute
\[\ll_\delta q^{\xi+\delta} e^{(2\delta+\delta^2)y}.\]
Finally the integral from $1-\theta-i$ to $1-\theta+i$ is
\[\ll_\delta q^{\xi+\delta}\exp\{(1-\theta)^2y\}\int_{-1}^1 e^{-t^2 y}dt
\ll_\delta q^{\xi+\delta} \exp\{(1-\theta)^2y\}.\]
Since $\theta<\tfrac12$ we therefore deduce that
\[I=\prod_{p|q}(1-p^{-1})L(1,\chi)f_0(1)e^y
+O_\delta\left(q^{\xi+\delta} \exp\{(1-\theta)^2y\}\right)\]
provided that $\delta<\tfrac{1}{10}$, say. Comparing the error term here 
with the lower bound (\ref{pb}) for the main term, we deduce that
\beql{IB2}
I\gg_\delta q^{\kappa-\delta},
\eeq
provided that
\beql{CC}
\delta\le\tfrac13 (\kappa(2\theta-\theta^2)-\xi),
\eeq
for example. We can always choose a suitable $\delta$ whenever
$\xi<\kappa(2\theta-\theta^2)$.

We now examine the infinite sum in (\ref{I2}) and show that terms with
$n$  much larger than $e^{2y}$ make a negligible contribution. For
the range $n> e^{8y}$ we note that $a_n\le n$, whence
\begin{eqnarray*}
\sum_{n>e^{8y}}a_n\exp\{-(\log n)^2/4y\} & \ll &
\int_{e^{8y}}^\infty x\exp\{-(\log x)^2/4y\}dx\\
&=& \int_{8y}^\infty \exp\{2u-u^2/4y\}du\\
&=& \int_{8y}^\infty \exp\{4y-(u-4y)^2/4y\}du.
\end{eqnarray*}
However $(u-4y)^2\ge 4y(u-4y)$ for $u\ge 8y$, so that
\[\sum_{n>e^{8y}}a_n\exp\{-(\log n)^2/4y\} \ll
e^{4y}\int_{8y}^\infty \exp\{-(u-4y)\}du=1.\]
This is negligible compared to $I$ by (\ref{IB2}), provided (\ref{CC}) holds.
We employ a similar argument for the range $e^{(2+\eta)y}<n\le e^{8y}$, 
but using the bound $a_n\ll_\delta n^\delta$. Here $\eta>0$ is a small 
constant to be chosen shortly. This time we find that
\begin{eqnarray*}
\sum_{e^{(2+\eta)y}<n\le e^{8y}}a_n\exp\{-(\log n)^2/4y\} & \ll_\delta &
e^{8y\delta}\int_{e^{(2+\eta)y}}^{e^{8y}} \exp\{-(\log x)^2/4y\}dx\\
&=& e^{8y\delta}\int_{(2+\eta)y}^{8y}\exp\{u-u^2/4y\}du\\
&=& e^{8y\delta}\int_{(2+\eta)y}^{8y}\exp\{y-(u-2y)^2/4y\}du\\
&\le & e^{8y\delta+y}\int_{(2+\eta)y}^{8y}\exp\{-\tfrac{\eta}{4}(u-2y)\}du\\
&\ll_{\delta,\eta}& e^{8y\delta+y-\eta^2 y/4}.
\end{eqnarray*}
In the same way
\begin{eqnarray*}
\sum_{n\le e^{(2-\eta)y}}a_n\exp\{-(\log n)^2/4y\} & \ll_\delta &
e^{8y\delta}\int_0^{e^{(2-\eta)y}} \exp\{-(\log x)^2/4y\}dx\\
&=& e^{8y\delta}\int_0^{(2-\eta)y}\exp\{u-u^2/4y\}du\\
&=& e^{8y\delta}\int_0^{(2-\eta)y}\exp\{y-(u-2y)^2/4y\}du\\
&\le & e^{8y\delta+y}\int_0^{(2-\eta)y}\exp\{-\tfrac{\eta}{4}|u-2y|\}du\\
&\ll_{\delta,\eta}& e^{8y\delta+y-\eta^2 y/4}.
\end{eqnarray*}
We now make the choice $\eta=\delta^{1/3}$. The bound (\ref{IB2}) then 
shows that the sums for $n>e^{(2+\eta)y}$ and $n\le e^{(2-\eta)y}$ are 
negligible compared to $I$
when $\delta<\kappa$ and $8\kappa\delta-\delta^{2/3}\kappa/4<-\delta$.
We therefore impose the additional condition that
\beql{CC2}
\delta\le \frac{\kappa^3}{65(1+8\kappa)^3}.
\eeq
Thus, under the conditions (\ref{CC}) and (\ref{CC2}) it follows that
\[\frac{1}{2\sqrt{\pi y}}\sum_{e^{(2-\eta)y}<n\le e^{(2+\eta)y}}
a_n\exp\{-(\log n)^2/4y\}
\gg_{\kappa,\theta,\xi} q^{\kappa-\delta},\]
whence
\[\sum_{n\le e^{(2+\eta)y}}a_n \gg_{\kappa,\theta,\xi} 
q^{\kappa-\delta}\exp\{(2-\eta)^2y/4\},\]

In order to apply Lemma \ref{FundLemV} we choose 
$\varpi=1/(2\ell)-\delta$ and $\kappa=\varpi/(2+\eta)$, where 
$\eta=\delta^{1/3}$. We then obtain 
\[M\gg_{\ell,\delta,\theta,\xi}  q^{\kappa(1+(2-\eta)^2/4)-\delta},\]
so that we obtain $h_\ell(K)\ll_{\ep} |\Delta_K|^{1/2-1/2\ell+\ep}$ 
by choosing $\delta$ sufficiently small.

\section{Cubic Fields --- Proof of Theorem \ref{CT}}

For our proof of Theorem \ref{CT} we write $q$ for the conductor of
$\chi$. We will introduce a constant $\delta$ in the range
$0<\delta<1/(4\ell)$. 
All the order constants in this section will be effective, and will be
allowed to
depend on $\ell$ and $\delta$, but any other dependencies will be explicitly
mentioned. We define $y=|\Delta_K|^{1/4\ell-\delta}$ and
\begin{eqnarray*}
f(s)&=&
(1-\delta)\int_{y^{1-\delta}}^y
\int_{z^{1-\delta}}^z u^{s-1}\frac{du}{u}\, \frac{dz}{z}\\
&=&
\frac{(1-\delta)y^{s-1}-(2-\delta)y^{(1-\delta)(s-1)}+y^{(1-\delta)^2(s-1)}}{(s-1)^2}.
\end{eqnarray*}
We see that
\[f(1)=\tfrac12\delta^2(1-\delta)(2-\delta),\]
that $f(s)\ge 0$ for all real $s$, and that $f(s)\ll |s-1|^{-2}$
when $\mathrm{Re}(s)\le 1$. 

If we now define
\[I=\frac{1}{2\pi i}\int_{2-i\infty}^{2+i\infty}
\left\{-\frac{L'(s,\chi)}{L(s,\chi)}\right\} f(s)ds\]
it follows from termwise integration that
\[I=\sum_{y^{(1-\delta)^2}<n\le y}\frac{\chi(n)\Lambda(n)}{n}
\min\left(\log(n/y^{(1-\delta)^2})\,,\,\log(y^{1-\delta}/n^{1-\delta})\right).\]
By the functional equation we have 
\[\frac{L'(-\tfrac12+it,\chi)}{L(-\tfrac12+it,\chi)}\ll \log q(|t|+2).\]
Hence on shifting the line of integration to $\sigma=-\tfrac12$ we
find that
\[ I=-\sum_\rho f(\rho)+
O\left(\int_{-\infty}^\infty \log q(|t|+2) \frac{dt}{1+t^2}\right)
=-\sum_\rho f(\rho)+O(\log q).\]
Zeros $\rho$ with $|\gamma|\ge 1$ contribute $O(\log q)$, whence
\[I=-\sum_{\rho,\, |\gamma|\le 1}f(\rho)+O(\log q).\]
A precisely analogous argument using $\zeta(s)$ in place of
$L(s,\chi)$ yields
\[\sum_{y^{(1-\delta)^2}<n\le y}\frac{\Lambda(n)}{n}
\min\left(\log(n/y^{(1-\delta)^2})\,,\,\log(y^{1-\delta}/n^{1-\delta})\right)
=f(1)(\log y)^2+O(1),\]
since $\zeta(s)$ has a simple pole at $s=1$, but no zeros with
$|\gamma|\le 1$. We may remove terms $n=p^e$ with $e\ge 2$ at a cost
$O(1)$; and also terms $n=p$ for which $p|q$, at a cost $O(\log q)$,
so that
\begin{eqnarray}\label{psl}
\lefteqn{\sum_{y^{(1-\delta)^2}<p\le y}\frac{\chi(p)\log p}{p}
\min\left(\log(p/y^{(1-\delta)^2})\,,\,\log(y^{1-\delta}/p^{1-\delta})\right)}
\hspace{5cm}\nonumber\\
&=&-\sum_\rho f(\rho)+O(\log q)
\end{eqnarray}
and
\begin{eqnarray*}
\lefteqn{\sum_{y^{(1-\delta)^2}<p\le y}\frac{\chi_0(p)\log p}{p}
\min\left(\log(p/y^{(1-\delta)^2})\,,\,\log(y^{1-\delta}/p^{1-\delta})\right)}
\hspace{5cm}\\
&=& f(0)(\log y)^2+O(\log q).
\end{eqnarray*}

We proceed to consider the contribution to (\ref{psl}) arising from zeros 
for which $|\gamma|\le 1$, and we begin by considering the case in 
which $\chi$ is cubic. The relevant region 
may be covered by annuli $R/\log q<|\rho-1|\le 2R/\log q$, where $R$ 
runs over powers of 2 such that
\beql{ran}
\frac{\vartheta}{\log|\Delta_K|}\ll \frac{R}{\log q}\ll 1.
\eeq
According to Linnik's ``Density Theorem" (see Prachar 
\cite[Page 331]{Prachar}, for example) each such annulus contains
$O(R+1)$ zeros, and each zero contributes $O(R^{-2}(\log q)^2)$ to
(\ref{psl}). If we let $R$ run over powers of 2 in the range (\ref{ran})
we therefore find that the overall contribution to (3.1) is
\[\ll \vartheta^{-2}(\log|\Delta_K|)^2+\vartheta^{-1}(\log|\Delta_K|)(\log q)
\ll\vartheta^{-1}(\log|\Delta_K|)^2.\]
It follows that
\[\sum_{y^{(1-\delta)^2}<p\le y}\frac{\chi(p)\log p}{p}
\min\left(\log(p/y^{(1-\delta)^2})\,,\,\log(y^{1-\delta}/p^{1-\delta})\right)
\ll \vartheta^{-1}(\log|\Delta_K|)^2.\]
By complex conjugation the same estimate holds when $\chi$ is 
replaced by $\overline{\chi}$, and we deduce that
\begin{eqnarray*}
\lefteqn{\sum_{y^{(1-\delta)^2}<p\le y}
\frac{(\chi_0(p)+\chi(p)+\overline{\chi}(p))\log p}{p}
\min\left(\log(p/y^{(1-\delta)^2})\,,\,\log(y^{1-\delta}/p^{1-\delta})\right)}
\hspace{5.5cm}\\
&=& f(1)(\log y)^2+O(\vartheta^{-1}(\log|\Delta_K|)^2).
\end{eqnarray*}
We recall now that $y=|\Delta_K|^{1/4\ell-\delta}$, and that our
implied constants are allowed to depend on $\ell$ and $\delta$, but 
not on anything else.
Thus, if $\vartheta$ is sufficiently large in terms of $\ell$ and $\delta$
we will have
\[3\sum_{\substack{y^{(1-\delta)^2}<p\le y\\ \chi(p)=1}}
\frac{\log p}{p}
\min\left(\log(p/y^{(1-\delta)^2})\,,\,\log(y^{1-\delta}/p^{1-\delta})\right)
\gg(\log y)^2,\]
and hence
\[\sum_{\substack{y^{(1-\delta)^2}<p\le y\\ \chi(p)=1}}p^{-1}\gg 1\]
and
\[\sum_{\substack{y^{(1-\delta)^2}<p\le y\\ \chi(p)=1}}1
\gg y^{(1-\delta)^2}.\]
Thus Lemma \ref{FundLem}, with $\ep=\delta$, produces
\[h_\ell(K)\ll_{\ell,\delta}|\Delta_K|^{1/2-(1-\delta)^2(1/4\ell-\delta)+\delta}.\]
Taking $\delta$ suitably small gives the required conclusion, when 
$\chi$ has order 3.

When $\chi$ is quadratic we follow the same argument, but bearing in 
mind that $L(s,\chi)$ may have real zeros in the disc 
$|s-1|\le\vartheta/\log|\Delta|_K$. What we find is that
\begin{eqnarray*}
\lefteqn{\sum_{y^{(1-\delta)^2}<p\le y}
\frac{(\chi_0(p)-\chi(p))\log p}{p}
\min\left(\log(p/y^{(1-\delta)^2})\,,\,\log(y^{1-\delta}/p^{1-\delta})\right)}
\hspace{4cm}\\
&=& f(0)(\log y)^2+\sum_{\rho\;\mathrm{real}}f(\rho)
+O(\vartheta^{-1}(\log|\Delta_K|)^2).
\end{eqnarray*}
Since $f(\rho)\ge 0$ for real $\rho$ we deduce as above that
\[\sum_{\substack{y^{(1-\delta)^2}<p\le y\\ \chi(p)=-1}}
\frac{\log p}{p}
\min\left(\log(p/y^{(1-\delta)^2})\,,\,\log(y^{1-\delta}/p^{1-\delta})\right)
\gg(\log y)^2,\]
and then that
\[\sum_{\substack{y^{(1-\delta)^2}<p\le y\\ \chi(p)=-1}}1
\gg y^{(1-\delta)^2}.\]
The argument may now be completed as before.

\section{Proof of Theorem \ref{zfr}}
  As explained in the introduction, the proof we shall give is a
  straightforward amalgam of familiar arguments. We begin with a result
  due to Landau \cite{Landau2}, see Titchmarsh \cite[Lemma $\beta$,
    page 56]{Titch}.
  \begin{lemma}\label{landau}
There is an effectively computable constant $A_1>0$ with the following
property. Suppose $f(s)$ is regular on the disc $|s-s_0|\leq r$, and
    satisfies
    \[\left|\frac{f(s)}{f(s_0)}\right|\le e^M\]
    there, for some $M>1$. Suppose further that $f(s)$ has no zeros in
    the right hand half of the disc $|s-s_0|\leq r$.  Then
\[-\mathrm{Re}\left\{\frac{f'(s_0)}{f(s_0)}\right\}\le
\frac{A_1 M}{r}-\sum_{|\rho-s_0|\le\tfrac12 r}
\mathrm{Re}\left\{\frac{1}{s_0-\rho}\right\},\]
where $\rho$ runs over zeros of $f$, counted according to
multiplicity.  In particular,
\[-\mathrm{Re}\left\{\frac{f'(s_0)}{f(s_0)}\right\}\le\frac{A_1 M}{r}.\]
  \end{lemma}
The form of the result is not precisely as stated by 
Titchmarsh, but it is given as \cite[top of page 57]{Titch}.

We begin our treatment of Theorem \ref{zfr} by considering the case in
which $\chi$ has order 2, and $|t|\le c\theta/\phi$ for a suitable
numerical constant $c$.  We will use parameters $\varphi$ and
$\sigma_0$ which we will  
specify in due course so as to satisfy the conditions
\beql{s0c}
\varphi\ge\phi+1\;\;\;\mbox{and}\;\;\; 1+e^{-\varphi}\le\sigma_0\le 1+\tfrac14\theta.
\eeq
We proceed to examine
\[-\frac{\zeta'(\sigma_0)}{\zeta(\sigma_0)}
-\frac{L'(\sigma_0,\chi)}{L(\sigma_0,\chi)}
=\sum_{n=1}^\infty\left(1+\chi(n)\right)\frac{\Lambda(n)}{n^{\sigma_0}}
\ge 0.\]
For the first term we have
$\zeta'(\sigma_0)/\zeta(\sigma_0)=-1/(\sigma_0-1)+O(1)$.
For the second term we apply Lemma \ref{landau} to the function
$L(s,\chi)$, taking
$s_0=\sigma_0$ and $r=\theta$, so that $L(s,\chi)$ has no zeros in
the right hand half of the disc $|s-s_0|\leq r$.  This disc is
contained in the region considered in (\ref{zfrc}), and
\[L(\sigma_0,\chi)^{-1}\le\zeta(\sigma_0)\le
(\sigma_0-1)^{-1}+1\le e^\phi +1\le e^{\varphi}.\]
so that 
\[\left|\frac{L(s,\chi)}{L(s_0,\chi)}\right|\le e^{\phi+\varphi}\le e^{2\varphi}.\]
We may therefore apply Lemma \ref{landau} with $M=2\varphi$, allowing us to
conclude that
\begin{eqnarray*}
0&\le&-\frac{\zeta'(\sigma_0)}{\zeta(\sigma_0)}
-\frac{L'(\sigma_0,\chi)}{L(\sigma_0,\chi)}\\
&\le& \frac{2A_1\varphi}{\theta}+O(1)+\frac{1}{\sigma_0-1}-
\sum_{|\rho-\sigma_0|\le\tfrac12\theta}
\mathrm{Re}\left\{\frac{1}{\sigma_0-\rho}\right\}.
\end{eqnarray*}
Since $\varphi\ge\phi\ge 1\ge \theta$ it follows that
\[0\le \frac{A_2\varphi}{\theta}+\frac{1}{\sigma_0-1}-
\sum_{|\rho-\sigma_0|\le\tfrac12\theta}
\mathrm{Re}\left\{\frac{1}{\sigma_0-\rho}\right\}\]
for a suitable constant $A_2\ge 1$.

However
\[\mathrm{Re}\left\{\frac{1}{\sigma_0-\rho}\right\}=
\frac{\sigma_0-\beta}{(\sigma_0-\beta)^2+\gamma^2}.\]
This is always non-negative. Moreover, if 
$|\gamma|\le\tfrac12 (\sigma_0-1)$ we have 
$|\gamma|\le\tfrac12 (\sigma_0-\beta)$ and hence
\[\mathrm{Re}\left\{\frac{1}{\sigma_0-\rho}\right\}\ge
\frac{\sigma_0-\beta}{(\sigma_0-\beta)^2+\tfrac14(\sigma_0-\beta)^2}
=\frac{4/5}{\sigma_0-\beta}.\]
It follows that
\[0\le \frac{A_2\varphi}{\theta}+\frac{1}{\sigma_0-1}-
\sum_{\substack{|\rho-\sigma_0|\le\tfrac12\theta\\ |\gamma|\le\tfrac12 (\sigma_0-1)}} \frac{4/5}{\sigma_0-\beta}.\]
If the sum includes two real zeros with 
$1-\beta<\tfrac12(\sigma_0-1)$, or a pair of complex conjugate
zeros in the above region, we deduce that
\[0<\frac{A_2\varphi}{\theta}+\frac{1}{\sigma_0-1}-
\frac{16/15}{\sigma_0-1},\]
and therefore that
\[\sigma_0-1>\frac{\theta}{15A_2\varphi}.\]
We can therefore obtain a contradiction by making the choice
\[\sigma_0=1+\frac{\theta}{15A_2\varphi}\]
provided that this satisfies the constraints (\ref{s0c}).  
Since $0<\theta\le 1\le\phi\le\varphi$ and $A_2\ge 1$ it is automatic that 
$\sigma_0\le 1+\tfrac14\theta$. Moreover, in view of our assumption that
$\theta\ge \phi e^{-\phi}$ we will have 
\[\sigma_0\ge 1+\frac{\phi e^{-\phi}}{15A_2\varphi}\ge 1+e^{-\varphi}\]
provided that we choose $\varphi=\phi+\kappa$ with a constant 
$\kappa\ge 1$ so large that
\[e^\kappa\ge 15A_2(1+\kappa)\ge 15A_2(1+\kappa/\phi).\]

We conclude from the above that there can be at most one zero 
lying in the disc $|\rho-\sigma_0|\le\tfrac12\theta$ such that
$1-\beta<\tfrac12(\sigma_0-1)$ and $|\gamma|\le\tfrac12(\sigma_0-1)$.
Any zero satisfying these latter two conditions will automatically have
\[|\rho-\sigma_0|\le \sigma_0-\beta+|\gamma|\le 2(\sigma_0-1)\le\tfrac12\theta,\]
in view of (\ref{s0c}). Since
\[\tfrac12(\sigma_0-1)=\frac{\theta}{30A_2\varphi}=
\frac{\theta}{30A_2(\phi+\kappa)}\ge\frac{\theta}{30A_2(\kappa+1)\phi}\]
this suffices for the theorem, with $C_2=(31A_2(\kappa+1))^{-1}$,
when $\chi$ is real and $|\gamma|\le C_2\theta/ \phi$.
\bigskip

The argument for the remaining cases of the theorem follows along very 
similar lines, and so we shall be brief.
We look next at the situation in which either $\chi$ is real with
$|\gamma|\ge C_2\theta/ \phi$, or $\chi$ is complex. Suppose
that $L(s,\chi)$ has a zero $\beta+i\gamma$ with $|\gamma|\le 1$.
We choose $\sigma_0$ and $\varphi$ as in (\ref{s0c}), and consider
the familiar inequality
\beql{3t}
-3\frac{\zeta'(\sigma_0)}{\zeta(\sigma_0)}
-4\mathrm{Re}
\left\{\frac{L'(\sigma_0+i\gamma,\chi)}{L(\sigma_0+i\gamma,\chi)}\right\}
-\mathrm{Re}
\left\{\frac{L'(\sigma_0+2i\gamma,\chi^2)}{L(\sigma_0+2i\gamma,\chi^2)}\right\}
\ge 0.
\eeq
We again apply Lemma \ref{landau} to the function $L(s,\chi)$, taking
$r=\theta$ as before, but this time choosing $s_0=\sigma_0+i\gamma$.
The previous argument then shows that
\begin{eqnarray*}
  -3\frac{\zeta'(\sigma_0)}{\zeta(\sigma_0)}
-4\mathrm{Re}
\left\{\frac{L'(\sigma_0+i\gamma,\chi)}{L(\sigma_0+i\gamma,\chi)}\right\}
&\le &\frac{8A_1\varphi}{\theta}+O(1)+\frac{3}{\sigma_0-1}\\
&&\hspace{1cm}-4\sum_{|\rho-s_0|\le\tfrac12\theta}
\mathrm{Re}\left\{\frac{1}{\sigma_0-\rho}\right\}.
\end{eqnarray*}
If $\rho=\beta+i\gamma$ is
included it contributes 
$(\sigma_0-\beta)^{-1}$ to the sum. Moreover any other zero will make a 
non-negative contribution, so that
\begin{eqnarray*}
  \lefteqn{-3\mathrm{Re}
    \left\{\frac{L'(\sigma_0,\chi_0)}{L(\sigma_0,\chi_0)}\right\}
    -4\mathrm{Re}
    \left\{\frac{L'(\sigma_0+i\gamma,\chi)}{L(\sigma_0+i\gamma,\chi)}\right\}}
  \hspace{2cm}\\
&\le& \frac{8A_1\varphi}{\theta}+O(1)+\frac{3}{\sigma_0-1}-
\frac{4}{\sigma_0-\beta},
\end{eqnarray*}
provided that $\sigma_0-\beta\le\tfrac12\theta$.

It remains to consider the third term in (\ref{3t}). When $\chi$ is
complex we may apply
Lemma \ref{landau} to $L(s,\chi^2)$, taking $s_0=\sigma_0+2i\gamma$
and $r=\theta$. This produces
\[-\mathrm{Re}
\left\{\frac{L'(\sigma_0+2i\gamma,\chi^2)}{L(\sigma_0+2i\gamma,\chi^2)}\right\}
\le\frac{2A_1\varphi}{\theta}.\]
Alternatively, if $\chi$ is real and $|\gamma|\ge C_2\theta/ \phi$,
then $\chi^2=\chi_0$ and 
\begin{eqnarray*}
  \left|\frac{L'(\sigma_0+2i\gamma,\chi^2)}{L(\sigma_0+2i\gamma,\chi^2)}\right|
&=& \left|\sum_1^\infty\chi_0(n)\Lambda(n)n^{-\sigma_0-2i\gamma}\right|\\
&\le&
\left|\frac{\zeta'(\sigma_0+2i\gamma)}{\zeta(\sigma_0+2i\gamma)}\right|+
\sum_{p|q}\sum_{k=1}^\infty \frac{\log p}{p^{k\sigma_0}}.
\end{eqnarray*}
Here we have
\[\left|\frac{\zeta'(\sigma_0+2i\gamma)}{\zeta(\sigma_0+2i\gamma)}\right|=
\frac{1}{|\sigma+2i\gamma-1|}+O(1)\le\frac{1}{2|\gamma|}+
O(1)\le \frac{\phi}{2C_2\theta}+O(1)\]
and
\beql{ll3}
\sum_{p|q}\sum_{k=1}^\infty \frac{\log p}{p^{k\sigma_0}}\le
\sum_{p|q}\frac{\log p}{p-1}\ll \log\log 3q.
\eeq
In this case, since $\phi\le\varphi$ we have
\[-\mathrm{Re}
\left\{\frac{L'(\sigma_0+2i\gamma,\chi^2)}{L(\sigma_0+2i\gamma,\chi^2)}\right\}
\le\frac{\varphi}{2C_2\theta}+O(\log\log 3q).\]
Thus, in either case we may conclude that
\[0\le\frac{A_3\varphi}{\theta}+\frac{3}{\sigma_0-1}-\frac{4}{\sigma_0-\beta}
+O(\log\log 3q)\]
for a suitable constant $A_3$, provided that $\sigma_0-\beta\le\tfrac12\theta$. 
We assume in Theorem \ref{zfr} that 
$\phi\ge\theta(1+\log\log 3q)$, whence
\[0\le\frac{A_4\varphi}{\theta}+\frac{3}{\sigma_0-1}-\frac{4}{\sigma_0-\beta}\]
with a further constant $A_4$. If $1-\beta\le\tfrac14(\sigma_0-1)$
we deduce that
\[0\le\frac{A_4\varphi}{\theta}-\frac{1}{15(\sigma_0-1)},\]
which gives the desired contradiction on choosing
\[\sigma_0=1+\frac{\theta}{16A_4\varphi}.\]
We can then complete the argument as before, taking $\varphi=\phi+\kappa$
with a suitable constant $\kappa\ge 1$, and checking that the assumption
$1-\beta\le\tfrac14(\sigma_0-1)$ ensures that we have $\sigma_0-\beta\le\tfrac12\theta$.
We conclude that, for the two cases under consideration, there are no zeros
with $\beta\ge 1-C_2'\theta/\phi$, for a new constant $C_2'>0$.

Finally we consider two different non-principal real characters
$\chi_1$ and $\chi_2$, using the inequality
\[-\frac{\zeta'(\sigma_0)}{\zeta(\sigma_0)}-
\frac{L'(\sigma_0,\chi_1)}{L(\sigma_0,\chi_1)}
-\frac{L'(\sigma_0,\chi_2)}{L(\sigma_0,\chi_2)}-
\frac{L'(\sigma_0,\chi_1\chi_2)}{L(\sigma_0,\chi_1\chi_2)}
\ge 0.\]
If $L(s,\chi_1)$ and $L(s,\chi_2)$ have real zeros at $\beta_1$
and $\beta_2$ respectively, 
both in the range $\sigma_0-\beta\le\tfrac12\theta$, then the
previous arguments lead to the bound
\[0\le\frac{6A_1\varphi}{\theta}+O(1)+\frac{1}{\sigma_0-1}-
\frac{1}{\sigma_0-\beta_1}
-\frac{1}{\sigma_0-\beta_2}.\]
If $1-\beta_1\le\tfrac12(\sigma_0-1)$ and
$1-\beta_2\le\tfrac12(\sigma_0-1)$ we then obtain a 
contradiction as before, on making a suitable choice of $\sigma_0$.
This shows that
there can be at most one real character $\chi$ with a zero in the
range $\beta\ge 1-C_2\theta/\phi$, thereby completing the proof of the theorem.

Mathematical Institute, 

Radcliffe Observatory Quarter, 

Woodstock Road, 

Oxford  OX2~6GG,

UK
\bigskip

 {\tt rhb@maths.ox.ac.uk}

\end{document}